\documentclass[12pt,a4paper,reqno]{article}
\usepackage{fleqn,amssymb,amscd,a4,enumerate,bm,amsthm,amsmath,float,multicol,}

\usepackage[a4paper,textwidth=5.25in,textheight=8.25in,centering]{geometry}
\theoremstyle{plain}
\newtheorem{theorem}{\bf Theorem}
\newtheorem*{Proposition}{\bf Proposition}

\theoremstyle{remark}

\def\RR{\mathbb{R}}
\usepackage{titlesec,soul}
\titleformat{\section}[block]{\large\bfseries}{\thesection.}{.5em}{}
\titlespacing*{\section}{0pt}{18pt}{15pt}
\date{}

\begin{document}

\thispagestyle{empty} \noindent {\large\bf Cohomology algebra of the
orbit space of free circle group actions on lens spaces}

\vskip1.5em
\noindent{\bf Hemant Kumar Singh and Tej Bahadur Singh}\\[.5em]
Department of Mathematics, University of Delhi, Delhi 110 007, India

\baselineskip24pt \vskip.5in\noindent \textbf{Abstract.} Suppose
that  $G =\mathbb{S}^1$ acts freely on a finitistic space $X$ whose
$\bmod\,\,p$ cohomology ring isomorphic to that of a lens space
$L^{2m-1}(p;q_1,\ldots,q_m)$. In this paper, we determine the
$\bmod\,\,p$ cohomology ring of the orbit space $X/G$. If the
characteristic class $\alpha \,\epsilon\, H^2(X/G;\mathbb{Z}_p)$ of
the $\mathbb{S}^1$-bundle $\mathbb{S}^1\hookrightarrow X\rightarrow
X/G$ is nonzero, then the $\bmod\,\,p$ index of the action is
defined to be the largest integer $n$ such that $\alpha^n \neq 0$.
We also show that the $\bmod\,\,p$ index of a free action of
$\mathbb{S}^1$ on a lens space $L^{2m-1}(p;q_1,\ldots,q_m)$ is
$p-1$, provided that $\alpha\neq 0$.

\bigskip
\bigskip
\bigskip\noindent
\textbf{Key Words:} Characteristic class, Finitistic space, Free
action, Index, Spectral sequence.

\bigskip\noindent
\textbf{2000 Mathematics Subject classification.} Primary
57S17; Secondary 57S25

\thispagestyle{empty}

\vfill
\eject
\baselineskip24pt

\section{Introduction}

Let $X$ be a topological space and $G$ be a topological group acting
continuously on $X$. The set $\hat{x} =\{gx|g \in G\}$ is called the
orbit of $x$. The set of all orbits $\hat{x}$, $x \in X$, is denoted
by $X/G$ and assigned the quotient topology induced by the natural
projection $\pi :X \to X/G$, $x \to \hat{x}$. An action of $G$ on
$X$ is said to be free if $g(x) =x$, for any $x \in X\Rightarrow g
=e$, the identity element of $G$. The orbit space of a free
transformation group $(G, S^n)$, where $G$ is a finite group, has
been studied extensively ([2], [7], [8], [10], [15]). However, a
little is known if the total space $X$ is a compact manifold other
than a sphere ([3], [6], [9], [14]). The orbit space of a free
involution on a real or complex projective space has been studied by
the authors in [13]. We have also determined the cohomology algebra
of the orbit space of free actions of $Z_p$ on a generalized lens
space $L^{2m-1}(p;q_1, q_2,\ldots, q_m)$ in [12]. In this note, we
determine the $\bmod\,\, p$ cohomology algebra of obrit spaces of
free actions of circle group $\mathbb{S}^1$ on the real projective
space and a lens space. Note that $\mathbb{S}^1$ can not act freely
on a 'finitistic' space having integral cohomology of a
finite-dimensional complex projective space or a quaternionic
projective space (Theorem 7.10 of Chapter III, [1]). We recall that
a paracompact Hausdorff space is finitistic if every open covering
has a finite-dimensional refinement.

Throughout this paper, $H^*(X)$ will denote the \v{C}ech  cohomology
of the space $X$. It is known that $H^*(\RR P^n;\mathbb{Z}_2)
=Z_2[a]/ \langle a^{n+1}\rangle$, where deg $a =1$, and
$H^*(L^{2m-1}(p;q_1,\ldots,q_m);\mathbb{Z}_p)=\wedge (a)\otimes
{\mathbb{Z}_p[b]}/{\langle b^m\rangle}$, $\deg a=1$, $\beta(a)=b$,
where $\beta$ is the Bockstein homomorphism associated with the
coefficient sequence
$0\rightarrow\mathbb{Z}_p\rightarrow\mathbb{Z}_{p^2}\rightarrow\mathbb{Z}_p\rightarrow
0$ . By $X\thicksim_p Y$, we shall mean that $H^*(X; \mathbb{Z}_p)$
and $H^*(Y; \mathbb{Z}_p)$ are isomorphic. We establish the
following results.

\begin{theorem} Let $G =\mathbb{S}^1$ act freely on a finitistic space $X \thicksim_ p L^{2m -1}$ $(p;q_1, q_2,\ldots, q_m),$ $p$ an odd
prime. Then $H^*(X/G;\mathbb{Z}_p)$ is one of the following graded
commutative algebras:
\begin{enumerate}[\rm(i)]
\item $\dfrac{\mathbb{Z}_p[x,y_1,y_3,\ldots, y_{2p-3}, z]}{\langle x^p, z^n, xy_q, y_q y_{q'} -A_{qq'} x^{\frac{q+q'}{2}} - B_{qq' } zx^{\frac{q+q' -2p}{2}} \rangle}$

where $m =np, \deg x =2$, $\deg y_q =q$, $\deg z =2p$, $A_{qq'} =0$
when $q+q' > 2p$, $B_{qq'} =0$ when $q+q' < 2p$ and both $A_{qq'}$
and $B_{qq'}$ are zero when $q =q'$ or $q + q' =2p$.

\item $\mathbb{Z}_p[z]/\langle z^m\rangle$, where $\deg z =2$.
\end{enumerate}\end{theorem}
For free actions  of circle group on a cohomology real projective space, we have

\begin{theorem} Let $G =\mathbb{S}^1$ act freely on a finitistic space $X \thicksim_2 \RR P^{2m-1}$. Then
\begin{eqnarray*}
H^*(X/G;\mathbb{Z}_2)\cong \mathbb{Z}_2[z]/\langle z^m\rangle \quad \text{where } \deg z =2.
\end{eqnarray*}
\end{theorem}
Let $G=\mathbb{S}^1$ act freely on a space $X$, then there is an
orientable 1-sphere bundle $\mathbb{S}^1\hookrightarrow
X\overset{\nu}\rightarrow X/G$, where $\nu$ denotes the orbit map.
Let $\alpha\,\epsilon\,H^2(X/G;\mathbb{Z})$ be its characteristic
class. Jaworowski [4] has defined the (integral) index of a free
$\mathbb{S}^1$-action on the space $X$ to be the largest integer $n$
(if it exists) such that $\alpha^n\neq0$. Similarly, one can define
mod $p$ index of a free $\mathbb{S}^1$-action on a space $X$.
Jaworowski has shown that the (integral or rational)
$\mathbb{S}^1$-index of $L^{2m -1}(p;q_1, q_2,\ldots, q_m)$ is
$m-1$. It follows from the Thom-Gysin sequence for bundle
$\mathbb{S}^1\hookrightarrow X\rightarrow X/G$ that the
characteristic class is nonzero only if $p>2$ on a space  $X=L^{2m
-1}(p;q_1, q_2,\ldots, q_m)$. In this case, the mod $p$
$\mathbb{S}^1$-index of $X$ is $p-1$. It should be noted that
$G=\mathbb{S}^1$ can not act freely on $X \thicksim_2 \RR P^{2m}$.

\section{Preliminaries}

Let $G=\mathbb{S}^1$ act on a paracompact Hausdorff space $X$. Then
there is an associated fibration $X \overset{i}{\longrightarrow} X_G
\overset{\pi}{\longrightarrow} B_G$, where $X_G =(E_G \times X)/G$
and $ E_G=\mathbb{S}^\infty\to B_G=\mathbb{C}P^\infty$ is a
universal $G$-bundle. It is known that $B_G$ is a CW-complex with
$2N$-skelton $\mathbb{C}P^N$ for all $N$ and $E_G$ is a CW-complex
with $2N+1$-skelton $\mathbb{S}^{2N+1}$. Write
$E^N_G=\mathbb{S}^{2N+1}$ and $B^N_G=\mathbb{C}P^N$. Then,
$H^i(E^N_G)=0$ for $0<i<2N+1$. Let $X^N_G=X\times_G E^N_G$ is
associated bundle over $B^N_G$ with fibre $X$. Then the equivariant
projection $X\times E^N_G\rightarrow X$ induces the map $\phi
:X_G^N\rightarrow X/G$. Let $G$ acts freely on $X$, then
\begin{eqnarray*}
\phi^*:H^i(X/G)\rightarrow H^i(X_G)
\end{eqnarray*}
is an isomorphism for all $i<2N+1$ with coefficient group
$\mathbb{Z}_p$, $p$ a prime, by Vietoris-Begle mapping theorem. By
$H^i(X_G)$ we mean $H^i(X^N_G)$, $N$ is large.

 To compute $H^*(X_G)$ we exploit the Leray-Serre spectral sequence of the map $\pi:X_G \to
B_G$ with coefficients in $\mathbb{Z}_p$, $p$ a prime. The edge
homomorphisms,
\begin{eqnarray*}
H^p(B_G) =E_2^{p,0} \to E_3^{p,0}\to \ldots \to E_{p+1}^{p,0}
=E_{\infty}^{p,0} \subseteq H^p (X_G),\text{and}
\end{eqnarray*}
\begin{eqnarray*}
H^q(X_G)\to E_{\infty}^{0,q} =E_{q+1}^{0,q}\subset \ldots \subset E_2^{0,q} = H^q(X)
\end{eqnarray*}
are the homomorphisms
\begin{eqnarray*}
\pi^*: H^p (B_G) \to H^p (X_G),\quad\text{and}\quad i^*: H^q (X_G)
\to H^q (X),
\end{eqnarray*}
respectively. We also recall the fact that the cup product in
$E_{r+1}$ is induced from that in $E_r$ via the isomorphism
$E_{r+1}\cong H^*(E_r)$. For the above facts, we refer to McCleary
[5].

\section{Proofs}
To prove our theorems, we need the following:
\begin{Proposition} Let $G =\mathbb{S}^1$ act freely on a finitistic space $X$ with $H^i(X)=0\, \forall i>n$.
Then $H^i(X/G)=0\, \forall i\geq n$ with coefficient group
$\mathbb{Z}_p$, $p$ a prime.
\end{Proposition}

 \proof[Proof:] We recall that the bundle $\mathbb{S}^1\hookrightarrow X\overset{\nu}\rightarrow
X/G$ is orientable, where $\nu:X\rightarrow X/G$ is the orbit map.
Consider, The Thom-Gysin sequence
\begin{eqnarray*}
\ldots\rightarrow
H^i(X/G)\overset{\nu^*}{\longrightarrow}H^i(X)\overset{\lambda^*}
{\longrightarrow}H^{i-1}(X/G)\overset{\mu^*}{\longrightarrow}H^{i+1}(X/G)\rightarrow
\ldots
\end{eqnarray*}
of the bundle, where $\mu^*$ is the multiplication by a
characteristic class $\alpha\,\epsilon\, H^2(X/G)$. This implies
that $H^i(X/G)\overset{\mu^*}{\longrightarrow}H^{i+2}(X/G)$ is an
isomorphism $\forall i\geq n$. Since $X$ is finitistic, $X/G$ is
finitistic (Deo and Tripathi [11]). Therefore, $H^*(X/G)$ can be
defined as the direct limit of $H^*(K(\cal{U}))$, where $K(\cal{U})$
denotes the nerve of $\cal{U}$ and $\cal{U}$ runs over all finite
dimensional open coverings of $X/G$. Let $\beta\,\epsilon\,H^i(X/G)$
be arbitrary. Then, we find a finite dimensional covering $\cal{V}$
of $X/G$ and elements $\alpha'\,\epsilon\,H^2(K(\cal{V}))$,
$\beta'\,\epsilon\,H^i(K(\cal{V}))$ such that $\rho(\alpha')=\alpha$
and $\rho(\beta')=\beta$ where $\rho
:\sum_{\cal{U}}H^i(K({\cal{U}})){\longrightarrow}H^i(X/G)$ is the
canonical map. Consequently, for $2k+i>$ dim $\cal{V}$, we have
$(\alpha')^k\beta'=0$ which implies that
$(\mu^*)^k(\beta)=\alpha^k\beta=0$. So, $\beta=0$, and the
proposition follows.

Now, we prove our main theorems.

 \proof[Proof of Theorem \emph{1}]
The case $m=1$ is trivial. So we assume $m>1$. Since
$G=\mathbb{S}^1$ acts freely on $X$, the Leray-Serre spectral
sequence of the map $\pi:X_G \to B_G$ does not collapse at the
$E_2$-term. As $\pi_1(B_G)$ is trivial, the fibration
$X\overset{i}{\rightarrow}X_G \overset{\pi}{\rightarrow}B_G$ has a
simple system of local coefficients on $B_G$. So the spectral
sequence has
\begin{eqnarray*}
E_2^{k,l } \cong H^k (B_G)\otimes H^{l } (X).
\end{eqnarray*}
Let $a \in H^1 (X)$ and $b\in H^2(X)$ be generators of the
cohomology ring $H^*(X)$ then $a^2 = 0$ and $b^m =0$. Consequently,
we have either ($d_2(1\otimes a) = t\otimes 1$ and $d_2(1\otimes b)
=0$) or $(d_2(1\otimes a) =0$ and $d_2(1\otimes b) =t\otimes a)$.
\noindent

\noindent
 {\textbf{Case I.}} If $d_2(1\otimes a) =0$ and
$d_2(1\otimes b) =t \otimes a$, then we have $d_2(1\otimes b^q) =qt
\otimes ab^{q-1}$ and $d_2(1\otimes ab^q) =0$ for $1\leq q < m$. So
$0 =d_2[(1\otimes b^{m-1})\cup (1\otimes b)]=mt\otimes ab^{m-1}$.
This forces $p|m$. Suppose that $m=np$. Now,
\begin{eqnarray*}
d_2:E_2^{k,l } \to E_2^{k+2,l  -1}
\end{eqnarray*}
is an isomorphism if $l $ is even and $2p$ does not divide $l $; and
the trivial homomorphism if $l $ is odd or $2p$ divides $l $. So
$E_3^{k,l } \cong E_2^{k,l } \cong\mathbb{Z}_p$ for even $k$ and $l
=2qp$ or $2(q+1)p-1$, $0\leq q < n$; $k=0, l $ is odd and $2p$ does
not divide $l $; and $E_3^{k,l} =0$, otherwise. Clearly, all the
differentials $d_3, d_4,\ldots, d_{2p-1}$ are trivial. Obviously,
\begin{eqnarray*}
d_{2p}: E_{2p}^{k,2qp}\to E_{2p}^{k+2p, 2(q-1)p+1}
\end{eqnarray*}
are the trivial homomorphisms for $q=1,2,\ldots, n-1$. If
\begin{eqnarray*}
d_{2p}: E_{2p}^{0,2p-1}\to E_{2p}^{2p,0}
\end{eqnarray*}
is also trivial, then
\begin{eqnarray*}
d_{2p}: E_{2p}^{k,2qp-1}\rightarrow E_{2p}^{k+2p,2(q-1)p}
\end{eqnarray*}
is the trivial homomorphism for $q =2,\ldots,n-1$, because every
element of $E_{2p}^{k,2qp-1}$ (even $k$) can be written as the
product of an element of $E_{2p}^{k,2(q-1)p}$ by $[1\otimes
ab^{p-1}] \in E_{2p}^{0,2p-1}$. It follows that $d_r = 0$, $\forall$
$r> 2p$ so that $E_{\infty} =E_3$. This contradicts the fact that
$H^i(X_G)=0$ for all $i\geq 2m-1$. Therefore,
\begin{eqnarray*}
d_{2p} : E_{2p}^{0,2p-1}\to E_{2p}^{2p,0}
\end{eqnarray*}
must be non-trivial. Assume that $d_{2p}([1\otimes ab^{p-1}]) =[t^p \otimes 1]$.
Then
\begin{eqnarray*}
d_{2p}:E_{2p}^{k,2qp-1}\to E_{2p}^{k+2p,2(q-1)p}
\end{eqnarray*}
is an isomorphism for all $k$ and $1 \leq q \leq n$. Now, it is
clear that $E_{\infty} =E_{2p+1}$. Also, $E_{2p+1}^{k,l}\cong
\mathbb{Z}_p$ for ((even) $k < 2p$, $l =2qp$, $(0\leq q < n)$) and
($k=0, l $ is odd and $2p$ does not divide $l$).  Thus
\begin{eqnarray*}
H^j(X_G)=\begin{cases} 0, & j =2qp -1(1\leq q \leq n) \ \text{or} \ \ j > 2np -2\\
\mathbb{Z}_p & \ \text{otherwise}. \end{cases}
\end{eqnarray*}
The elements $1\otimes b^p \in E_2^{0,2p}$ and $1 \otimes
ab^{(h-1)/2} \in E_2^{0,h}$, for $h =1,3,\ldots, 2p -3$ are
permanent cocycles. So they determine $z \in E_{\infty}^{0,2p}$ and
$y_q \in E_{\infty}^{0,q}$, $q =1,3,\ldots, 2p-3$, respectively.
Obviously, $i^*(z) =b^p$, $z^n =0$ and $y_q y_{q'} =0$. Let $x
=\pi^*(t) \in E_{\infty}^{2,0}$. Then $x^p =0$. It follows that the
total complex Tot $E_{\infty}^{*,*}$ is the graded commutative
algebra
\begin{eqnarray*}
Tot E_{\infty}^{*,*} =\dfrac{\mathbb{Z}_p[x,y_1, y_3,\ldots,
y_{2p-3},z]}{\langle x^p, y_q y_{q'}, x y_q, z^n\rangle}
\end{eqnarray*}
where $q$, $q' =1, 3,\ldots, 2p-3$.

Then $i^*(y_q) =ab^{\frac{(q-1)}{2}}$, $y_q^2 =0$ and $y_q y_{2p-q}
=0$. It follows that
\begin{eqnarray*}
H^*(X_G) =\dfrac{\mathbb{Z}_p[x, y_1, y_3,\ldots, y_{2p-3}, z]}
{\langle x^p,z^n, x y_q, y_q y_{q'} - A_{qq'} x^{\frac{q+q'}{2}}
-B_{qq'} z x^{\frac{q+q' -2p}{2}}\rangle}
\end{eqnarray*}
where $m =np$, $A_{qq'} = 0$ when $q + q' > 2p$, $B_{qq'} =0$ when
$q + q'< 2p$ and both $A_{qq'}$ and $ B_{qq'}$ are zero when $q =q'$
or $q + q' =2p$, $\deg x =2$, $\deg z =2p$, $\deg y_{q} =q$.

\noindent
 {\textbf{Case II.}} If $d_2(1\otimes a) =t\otimes 1 $ and
$d_2(1\otimes b) =0$, then
\begin{eqnarray*}
d_2: E_2^{k,l } \to E_2^{k+2,l  -1}
\end{eqnarray*}
is an isomorphism for $k$ even and $l$ odd and the trivial homomorphism for remaining values of $k$ and
$l $. Obviously, $E_3^{k,l } \cong \mathbb{Z}_p$ for $k=0$ and $l  =0,2,4,\ldots, 2 m-2$.
So that $E_{\infty} =E_3$. Therefore, we have
\begin{eqnarray*}
E_{\infty}^{k,l } =\begin{cases} \mathbb{Z}_p, & k=0 \ \text{and} \ \ l  =0,2,4,\ldots, 2m-2\\
0 & \ \text{otherwise}.\end{cases}
\end{eqnarray*}
The element $1 \otimes b \in E_2^{0,2}$ is a permanent cocycle and
determines an element $z \in E_{\infty}^{0,2}$. We have $i^*(z) =b$
and $z^m =0$. Therefore, the total complex Tot $E_{\infty}^{*,*}$ is
the graded commutative algebra
\begin{eqnarray*}
Tot E_{\infty}^{*,*} =\mathbb{Z}_p[z]/\langle z^m\rangle,  \quad
\deg z =2.
\end{eqnarray*}
It shows that $E_{\infty}^{0,l } = H^{l }(X_G)$ $\forall$ $l $ and
hence
\begin{eqnarray*}
H^*(X_G) =\mathbb{Z}_p[z]/\langle z^m\rangle, \quad \deg z =2.
\end{eqnarray*}
Since the action of $G$ on $X$ is free, the mod $p$ cohomology rings
of $X_G$ and $X/G$ are isomorphic. This completes the proof.

\noindent \proof[Proof of Theorem \emph{2}] As above we have
\begin{eqnarray*}
E_2^{kl } \cong H^k (B_G) \otimes H^{l } (X).
\end{eqnarray*}
Let $a \in H^1 (X)$ be the generator of the cohomology ring
$H^*(X)$. If $d_2(1\otimes a) =0$, then $d_2(1\otimes a^q) =0$, by
the multiplicative structure of spectral sequence. Consequently, the
spectral sequence degenerates which contradicts our hypothesis.
Therefore, we must have $d_2(1\otimes a) =t \otimes 1$. It is easily
seen that
\begin{eqnarray*}
d_2 : E_2^{k,l } \to E_2^{k+2, l-1}
\end{eqnarray*}
is an isomorphism for $k$ even and $l $ odd; and the trivial
homomorphism otherwise. So
\begin{eqnarray*}
E_{\infty}^{k,l }\cong \begin{cases} \mathbb{Z}_2, & k =0 \ \text{and} \ \ l  =0,2,4,\ldots, 2m -2\\
0,  & \text{otherwise} \end{cases}
\end{eqnarray*}
It follows that $H^*(X_G)$ and Tot $E_{\infty}^{*,*}$ are the same as
the graded commutative algebra. The case $m=1$ is
obvious, so assume that $m>1$.

The element $1\otimes a^2 \in E_2^{0,2}$ is a permanent cocycle and
determines an element $z \in E_{\infty}^{0,2} = H^2(X_G)$. We have
$i^*(z) =a^2$ and $z^m =0$. Therefore, the total complex Tot
$E_{\infty}^{*,*}$ is the graded commutative algebra.
\begin{eqnarray*}
Tot E_{\infty}^{*,*} =\mathbb{Z}_2[z]/\langle z^m\rangle, \ \text{where} \ \ \deg z=2.
\end{eqnarray*}
Thus $H^*(X_G) =\mathbb{Z}_2[z]/\langle z^m\rangle$, where $\deg z =2$. This
completes the proof.

\section{Examples}

Consider the $(2m-1)$ sphere $\mathbb{S}^{2m-1}\subset
\mathbb{C}\times \ldots\times \mathbb{C}$ $(m$ times). The map
$(\xi_1,\ldots,\xi_m)\to (z\xi_1,\ldots,z\xi_m)$, where $z\in
\mathbb{S}^1$, defines a free action of $G=\mathbb{S}^1$ on
$\mathbb{S}^{2m-1}$ with the orbit space
$\mathbb{S}^{2m-1}/\mathbb{S}^1$ the complex projective space. Let
$N=\langle z\rangle$, where $z=e^{2\pi i/p}$, then the orbit space
$\mathbb{S}^{2m-1}/N$ is the lens space $L^{2m-1}(p;1,\ldots,1)$
(resp. real projective space $\mathbb{R}P^{2m-1}$ for $p=2$). It
follows that there is a free action of $\mathbb{S}^1=G/N$ on a lens
space with the complex projective space as the orbit space. This
realizes the second case of Theorem $1$ and Theorem $2$. We note
that characteristic class of the bundle $\mathbb{S}^1\hookrightarrow
L^{2m -1}(p;q_1, q_2,\ldots, q_m)\rightarrow L^{2m -1}(p;q_1,
q_2,\ldots, q_m)/G$ over $\mathbb{Z}_p$ is zero. So, mod $p$ index
of this action is not defined.

\end{document}